\documentclass{article}
\title
{Multidimensional cellular automata and generalization of Fekete's lemma}
\author{Silvio Capobianco\footnote{
School of Computer Science, Reykjav\'{\i}k University;
email: \texttt{silvio@ru.is}
\newline
The author was partly supported by the project
``The Equational Logic of Parallel Processes''
(nr.~060013021) of The Icelandic Research Fund.}}
\date{}

\usepackage{amsfonts}

\newcommand{\integers}{\ensuremath{\mathbb{Z}}}
\newcommand{\naturals}{\ensuremath{\mathbb{N}}}
\newcommand{\reals}{\ensuremath{\mathbb{R}}}
\newcommand{\Acal}{\ensuremath{\mathcal{A}}}
\newcommand{\eg}{\textit{e.g.}}
\newcommand{\ie}{\textit{i.e.}}

\newcommand{\Neigh}{\ensuremath{\mathcal{N}}}
\newcommand{\Nset}{\ensuremath{\naturals}}
\newcommand{\Out}{\ensuremath{\mathrm{Out}}}
\newcommand{\Rset}{\ensuremath{\reals}}
\newcommand{\Xcal}{\ensuremath{\mathcal{X}}}
\newcommand{\Zcal}{\ensuremath{\mathcal{Z}}}
\newcommand{\Zset}{\ensuremath{\integers}}

\newtheorem{theorem}{Theorem}

\newenvironment{proof}{\textit{Proof.}}{$\Box$}

\begin{document}

\maketitle

\begin{abstract}
Fekete's lemma is a well known combinatorial result on number sequences:
we extend it to functions defined on $d$-tuples of integers.
As an application of the new variant,
we show that nonsurjective $d$-dimensional cellular automata
are characterized by loss of arbitrarily much information
on finite supports,
at a growth rate greater than that of the support's boundary
determined by the automaton's neighbourhood index.


\end{abstract}

\section{Introduction}

Let $f:\{1,2,\ldots\}\to[0,+\infty)$.
\textbf{Fekete's lemma}~\cite{fe23,lm95,vlw92} states that,
if $f(n+k)\leq f(n)+f(k)$ for every $n$ and $k$, then
\begin{equation} \label{eq_fekete-1}
\lim_{n\to\infty}\frac{f(n)}{n}=\inf_{n\geq 1}\frac{f(n)}{n}\;.
\end{equation}
The consequences of this simple statement are many and deep,
such as the definition of \emph{topological entropy}
for dynamical systems~\cite{akm65,ko01,lm95}
and Arratia's bound on the number of permutations
avoiding a given pattern~\cite{a99}.

More recently, in a joint work
with Tommaso Toffoli and Patrizia Mentrasti~\cite{tcm07},
we have made use of (\ref{eq_fekete-1}) to prove a result
on unidimensional, nonsurjective \emph{cellular automata} (CA).
CA are presentations of global dynamics in local terms:
each global state is a $d$-dimensional \emph{configuration},
and the global evolution rule changes the state locally at a site
by considering only the states of \emph{neighbouring} sites.
Nonsurjective CA are characterized by the existence
of a configuration \emph{on a finite region of the space}
that has no predecessor according to the evolution rule;
this can be restated from another point of view,
by saying that nonsurjective CA
\emph{lose information within finite range}.
Fekete's lemma then told us
that such CA must lose, on finite regions large enough,
an amount of information essentially proportional
to the size of the support itself---thus,
at least the size of the support's \emph{boundary}
determined by the neighborhood index;
since loss of information gives rise to \emph{extra channel capacity},
our team has devised a general algorithm to translate
from a presentation using an $n$-inputs, $1$-output local map
(\ie, CA) to one employing $n$-inputs, $n$-outputs events,
characteristic of a different class of presentations,
specifically, that of \emph{lattice gases} (LG).

In this paper, we state and prove a multivariate version of Fekete's lemma.
The motivation for this, is to provide a support to the conjecture
that the translation algorithm in~\cite{tcm07}
could be extended to arbitrary dimension.
To prove our generalization, we rearrange a proof of (\ref{eq_fekete-1})
so that it works on sequences of integer $d$-tuples,
after a suitable ordering on these is defined.
After that, we use the more general result at our hands to show
that the same phenomenon that allows rewriting CA as LG in dimension 1,
actually occurs in arbitrary finite dimension.
Incidentally, we get a criterion for CA surjectivity.

\section{Fekete's lemma, multivariate}

Let $\Zset_+=\Zset\cap(0,+\infty)$.
Consider the \textbf{product ordering} on $\Zset_+^d$
defined by $x\leq_\pi y$ iff $x_i\leq y_i$ for every $i\in\{1,\ldots,d\}$:
this is the kind of ordering used, \eg, in \emph{linear programming},
by writing $Ax\leq b$ to indicate a set of constraints
\begin{math}
a_{1,1}x_1+\ldots+a_{1,n}x_n\leq b_1,
\ldots,a_{m,1}x_1+\ldots+a_{m,n}x_n\leq b_m;
\end{math}
it is also the \emph{finest} ordering
that makes the projections \emph{monotonic}.
Observe that $\Zcal^d=(\Zset_+^d,\leq_\pi)$ is a \textbf{directed set},
\ie, for any two $x,y\in\Zset_+$ there exists $z\in\Zset_+$
such that both $x\leq_\pi z$ and $y\leq_\pi z$.
If $\Xcal=(X,\leq)$ is a directed set and $f:X\to\Rset$,
the \textbf{lower} and \textbf{upper limit} of $f$ in $\Xcal$
are defined as usual, \ie,
\begin{displaymath}
\liminf_{x\in\Xcal}f(x)
=\sup_{x\in X}\inf_{y\geq x}f(y)
\;\;\mathrm{and}\;\;
\limsup_{x\in\Xcal}f(x)
=\inf_{x\in X}\sup_{y\geq x}f(y)\;;
\end{displaymath}
moreover, $f$ has \textbf{limit} $L\in\Rset$ in $\Xcal$,
written $\lim_{x\in\Xcal}f(x)=L$,
if for every $\varepsilon>0$ there exists $x_\varepsilon\in X$
such that $|f(x)-L|<\varepsilon$ for every $x\geq x_\varepsilon$.
For example, if 
\begin{math}
r_1,\ldots,r_d\in\Nset=\Zset_+\cup\{0\}
\end{math}
are fixed, then
\begin{equation} \label{eq_ex-limit}
\lim_{(x_1,\ldots,x_d)\in\Zcal^d}
\frac{(x_1+r_1)\cdots(x_d+r_d)}{x_1\cdots x_d}=1\;.
\end{equation}
It follows from the definitions that
\begin{math}
\liminf_{x\in\Xcal}f(x)\leq\limsup_{x\in\Xcal}f(x)
\end{math},
and that $\lim_{x\in\Xcal}f(x)=L$ iff
\begin{math}
\liminf_{x\in\Xcal}f(x)=\limsup_{x\in\Xcal}f(x)=L
\end{math}.
\begin{theorem} \label{thm_fekete-d}
Let $f:\Zset_+^d\to[0,+\infty)$ satisfy
\begin{equation} \label{eq_subadd}
f(x_1,\ldots,x_j+y_j,\ldots,x_d)
\leq f(x_1,\ldots,x_j,\ldots,x_d)
+f(x_1,\ldots,y_j,\ldots,x_d)
\end{equation}
for every $x_1,\ldots,x_n,y_j\in\Zset_+$, $j\in\{1,\ldots,d\}$.
Then
\begin{equation} \label{eq_fekete-d}
\lim_{(x_1,\ldots,x_d)\in\Zcal^d}
\frac{f(x_1,\ldots,x_d)}{x_1\cdots x_d}
\end{equation}
exists, and equals
\begin{displaymath}
\inf_{x_1,\ldots,x_d\in\Zset_+}
\frac{f(x_1,\ldots,x_d)}{x_1\cdots x_d}\;.
\end{displaymath}
\end{theorem}
\begin{proof}
Because of (\ref{eq_subadd}), for every
\begin{math}
j\in\{1,\ldots,d\},\;x_1,\ldots,x_d\in\Zset_+,
\end{math}
if $x_j=qt+r$ with $q\in\Nset$ and $r\in\Zset_+$, then
\begin{equation} \label{eq_subadd-div}
f(x_1,\ldots,x_j,\ldots,x_d)
\leq qf(x_1,\ldots,t,\ldots,x_d)
+f(x_1,\ldots,r,\ldots,x_d)\;.
\end{equation}
Fix
\begin{math}
t_1,\ldots,t_d\in\Zset_+.
\end{math}
For each
\begin{math}
(x_1,\ldots,x_d)\in\Zset_+^d,
\end{math}
$d$ pairs
\begin{math}
(q_j,r_j)\in\Nset\times\Zset_+
\end{math}
are uniquely determined by
\begin{math}
x_j=q_jt_j+r_j
\end{math}
and
\begin{math}
1\leq r_j\leq t_j.
\end{math}
By repeatedly applying (\ref{eq_subadd-div}) to all of the $x_j$'s we find
\begin{equation} \label{eq_estimate}
\begin{array}{rcl}
f(x_1,\ldots,x_d) & \leq & q_1\cdots q_d f(t_1,\ldots,t_d) \\
& & + q_1\cdots q_{d-1} f(t_1,\ldots,t_{d-1},r_d) + \ldots \\
& & + q_1\cdots q_{d-2} f(t_1,\ldots,t_{d-2},r_{d-1},r_d) + \ldots \\
& & + \ldots
\end{array}
\end{equation}
where, in the next $k$'th line, $k\geq 1$,
each occurrence of $f$ has $k$ arguments chosen from the $r$'s
and $d-k$ chosen from the $t$'s,
and is multiplied precisely by the $q$'s corresponding to the $t$'s;
moreover, all these occurrences are bounded from above by the constant
\begin{math}
M=t_1\cdots t_d\cdot f(1,\ldots,1).
\end{math}
Divide both sides of (\ref{eq_estimate}) by $x_1\cdots x_d$:
since
\begin{math}
\lim_{x_j\to\infty}q_j/x_j=1/t_j,
\end{math}
if all the $x_j$'s are large enough,
then the first summand of right-hand side becomes very close to
\begin{math}
f(t_1,\ldots,t_d)/t_1\cdots t_d,
\end{math}
and the other ones become very small;
that is, for every $\varepsilon>0$, there exists
\begin{math}
(x_1,\ldots,x_d)\in X
\end{math}
such that, for every $(y_1,\ldots,y_d)\geq_\pi(x_1,\ldots,x_d)$,
\begin{displaymath}
\frac{f(y_1,\ldots,y_d)}{y_1\cdots y_d}
<\frac{f(t_1,\ldots,t_d)}{t_1\cdots t_d}+\varepsilon\;.
\end{displaymath}
From this follows
\begin{displaymath}
\limsup_{(x_1,\ldots,x_d)\in\Zcal^d}
\frac{f(x_1,\ldots,x_d)}{x_1\cdots x_d}
\leq\frac{f(t_1,\ldots,t_d)}{t_1\cdots t_d}\;;
\end{displaymath}
this is true whatever the $t_j$'s are, hence
\begin{displaymath}
\limsup_{(x_1,\ldots,x_d)\in\Zcal^d}
\frac{f(x_1,\ldots,x_d)}{x_1\cdots x_d}
\leq\inf_{t_1,\ldots,t_d\in\Zset_+}
\frac{f(t_1,\ldots,t_d)}{t_1\cdots t_d}\;.
\end{displaymath}
The thesis then follows from the inequality
\begin{displaymath}
\inf_{t_1,\ldots,t_d\in\Zset_+}
\frac{f(t_1,\ldots,t_d)}{t_1\cdots t_d}
\leq\liminf_{(x_1,\ldots,x_d)\in\Zcal^d}
\frac{f(x_1,\ldots,x_d)}{x_1\cdots x_d}\;.
\end{displaymath}
\end{proof}

\section{An application to cellular automata}

A \textbf{cellular automaton} (briefly, CA) is a quadruple
\begin{math}
\Acal=\left<d,Q,\Neigh,f\right>
\end{math}
where the \textbf{dimension} $d>0$ is an integer,
the \textbf{set of states} $Q$ is finite
and has at least two distinct elements,
the \textbf{neighbourhood index}
\begin{math}
\Neigh=\{\nu_1,\ldots,\nu_n\}
\end{math}
is a finite subset of $\Zset^d$,
and the \textbf{local evolution function} $f$ maps $Q^n$ into $Q$.
A \textbf{global evolution function} $F$ is induced by $f$
of the space $Q^{\Zset^d}$ of $d$-dimensional \textbf{configurations} by
\begin{equation} \label{eq:glob-evo}
(F(c))(x)=f\left(c(x+\nu_1),\ldots,c(x+\nu_n)\right)\;.
\end{equation}
$\Acal$ is said to be surjective if $F$ is.
For example, if $d=1$, $Q=\{0,1\}$, $\Neigh=\{+1\}$, $f(x)=x$, then
\begin{math}
\left<d,Q,\Neigh,f\right>
\end{math}
is the \emph{shift cellular automaton}
and
\begin{math}
(F(c))(x)=c(x+1)
\end{math}
is the \emph{shift map}, which is surjective;
on the other hand,
for same $d$ and $Q$,
\begin{math}
\Neigh=\{0,+1\},
\end{math}
and
\begin{math}
f(a,b)=a\cdot b,
\end{math}
we get a nonsurjective CA,
because if $\overline{c}(x)$ is 0 for $x=0$ and 1 otherwise, then
\begin{math}
F(c)\neq\overline{c}
\end{math}
for any $c$.

For every finite $E\subseteq\Zset^d$, calling
\begin{math}
E+\Neigh=\{x+\nu\mid x\in E,\nu\in\Neigh\},
\end{math}
a function
\begin{math}
F_E:Q^{E+\Neigh}\to Q^E
\end{math}
is induced by $f$, again by applying (\ref{eq:glob-evo}).
Observe that the number $|F_E(Q^{E+\Neigh})|$
of \textbf{patterns} over $E$ obtainable by applying (\ref{eq:glob-evo})
does not depend on the \emph{displacement} of $E$ along $\Zset^d$,
\ie, if $x+E=\{x+y\mid y\in E\}$,
then
\begin{math}
|F_{x+E}(Q^{x+E+\Neigh})|=|F_E(Q^{E+\Neigh})|.
\end{math}
It is well known (cf.~\cite{cal02})
that $\Acal$ is surjective
iff $F_E$ is surjective for every $E$
which is a \textbf{right $d$-polytope},
\ie, a subset of $\Zset^d$ of the form
\begin{math}
\prod_{i=1}^d\{k_i,\ldots,k_i+s_i-1\},
\end{math}
$s_1,\ldots,s_d\in\Zset_+$ being the \textbf{sides}.
(Here, ``right'' has the same meaning as in ``right triangle''.)
Put
\begin{math}
E(x_1,\ldots,x_d)=\prod_{i=1}^d\{0,\ldots,x_i-1\}:
\end{math}
if $\Neigh$ is contained in a right $d$-polytope of sides
\begin{math}
r_1,\ldots,r_d,
\end{math}
then
\begin{math}
E(x_1,\ldots,x_d)+\Neigh
\end{math}
is contained in a right $d$-polytope of sides
\begin{math}
x_1+r_1,\ldots,x_d+r_d,
\end{math}
which is the disjoint union of $E(x_1,\ldots,x_d)$ and a \textbf{boundary}.

Let
\begin{math}
\Acal=\left<d,Q,\Neigh,f\right>
\end{math}
be a CA.
If $\Acal$ is nonsurjective, then there must exist a support of suitable size
where not every possible pattern is reachable,
\ie, a part of the information is lost.
In the 1D case~\cite{tcm07}, such lost information
is proved to ultimately be as much as the boundary can transport;
which allowed devising a CA-to-LG conversion algorithm.
If the technique employed there is to be extended to higher dimension,
then we must determine whether such large a loss can still be achieved.

Call \textbf{output size} of $f$
on a right $d$-polytope of sides $x_1,\ldots,x_d$
the quantity
\begin{displaymath} 
\Out_f(x_1,\ldots,x_d)=
\left|F_{E(x_1,\ldots,x_d)}\left(Q^{E(x_1,\ldots,x_d)+\Neigh}\right)\right|\;.
\end{displaymath}
Then $\Acal$ is surjective iff
\begin{math}
\Out_f(x_1,\ldots,x_d)=|Q|^{x_1\cdots x_d}
\end{math}
for every $x_1,\ldots,x_d\in\Zset_+$.
By switching to a logarithmic measure unit,
we can associate to $\Acal$ a \textbf{loss of information}
on a right $d$-polytope of sides $x_1,\ldots,x_d$ defined as
\begin{equation} \label{eq_loi}
\Lambda_\Acal(x_1,\ldots,x_d)
=x_1\cdots x_d-\log_{|Q|}\Out_f(x_1,\ldots,x_d)\;.
\end{equation}
Observe how such loss is measured in \emph{$q$its} (with $q=|Q|$),
a $q$it being the amount of information carried by a $q$-states device;
$n$ $q$its correspond to $n\log_2 q$ bits.
\begin{theorem} \label{thm_surj}
Let $\Acal=\left<d,Q,\Neigh,f\right>$ be a CA.
Define $\Lambda_\Acal$ as by (\ref{eq_loi}).
Then
\begin{enumerate}
\item either $\Acal$ is surjective and $\Lambda_\Acal$ is identically zero,
\item or $\Acal$ is nonsurjective and 
for every
\begin{math}
K\geq 0,\;r_1,\ldots,r_d\in\Nset,
\end{math}
there exist $t_1,\ldots,t_d\in\Zset_+$
such that, for every
\begin{math}
(x_1,\ldots,x_d)\geq_\pi(t_1,\ldots,t_d),
\end{math}
\begin{displaymath}
\Lambda_\Acal(x_1,\ldots,x_d)
\geq(x_1+r_1)\cdots(x_d+r_d)-x_1\ldots x_d+K\;.
\end{displaymath}
\end{enumerate}
In particular, if $\Lambda_\Acal$ is bounded, then $\Acal$ is surjective.
\end{theorem}
\begin{proof}
Put $q=|Q|$.
Since a pattern over $E(x_1,\ldots,x_j+y_j,\ldots,x_d)$
can always be seen as the \emph{joining}
of a pattern over
\begin{math}
E(x_1,\ldots,x_j,\ldots,x_d)
\end{math}
and another one over
\begin{math}
E(x_1,\ldots,y_j,\ldots,x_d),
\end{math}
there cannot be more patterns obtainable over the former
than pairs of patterns obtainable over the latter,
\ie,
\begin{displaymath}
\Out_f(x_1,\ldots,x_j+y_j,\ldots,x_d)
\leq\Out_f(x_1,\ldots,x_j,\ldots,x_d)
\cdot\Out_f(x_1,\ldots,y_j,\ldots,x_d)
\end{displaymath}
whatever
\begin{math}
x_1,\ldots,x_n,y_j\in\Zset_+,\;j\in\{1,\ldots,d\}
\end{math}
are; consequently, $\log_q\Out_f$
is subadditive in each of its arguments (and nonnegative).
Let
\begin{equation} \label{eq_etaF}
\lambda_f
=\lim_{(x_1,\ldots,x_d)\in\Zcal^d}
\frac{\log_q\Out_f(x_1,\ldots,x_d)}{x_1\cdots x_d}\;,
\end{equation}
whose existence and value are given by Theorem~\ref{thm_fekete-d};
observe that 
$\lambda_f\leq 1$, and $\Acal$ is surjective iff $\lambda_f=1$.
Suppose $\Acal$ is nonsurjective.
Let $\delta\in(\lambda_f,1)$.
Choose $t_1,\ldots,t_d\in\Zset_+$ so that,
for every $(x_1,\ldots,x_d)\geq_\pi(t_1,\ldots,t_d)$, both
\begin{equation} \label{eq_variety}
\frac{\log_q\Out_f(x_1,\ldots,x_d)}{x_1\cdots x_d}
\leq\delta
\end{equation}
and
\begin{equation} \label{eq_boundary}
\frac{(x_1+r_1)\cdots(x_d+r_d)-x_1\cdots x_d+K}{x_1\cdots x_d}
\leq 1-\delta
\end{equation}
are satisfied,
the latter following from (\ref{eq_ex-limit}).
Then, for such $x_1,\ldots,x_d$,
\begin{eqnarray*}
x_1\cdots x_d-\log_q\Out_f(x_1,\ldots,x_d)
& \geq & (x_1\cdots x_d)(1-\delta) \\
& \geq & (x_1+r_1)\cdots(x_d+r_d)-x_1\cdots x_d+K\;.
\end{eqnarray*}
\end{proof}

From Theorem~\ref{thm_surj} follows that, for $(x_1,\ldots,x_d)$ satisfying
both (\ref{eq_variety}) and (\ref{eq_boundary}),
\emph{the loss of information is at least the size of the boundary}:
this is precisely the fact used in~\cite{tcm07},
and supports the conjecture that a similar construction
can be carried out in dimension $d>1$.
On the other hand---and perhaps, unfortunately---since
surjectivity of $d$-dimensional CA
is only decidable when $d=1$~\cite{ap72,ka90},
no algorithm exists to determine,
given an arbitrary multidimensional CA,
that its loss of information (\ref{eq_loi}) is bounded.

\section*{Acknowledgements}
The proof of Theorem~\ref{thm_fekete-d} is an adaptation of an argument
shown to us by Tullio Ceccherini--Silberstein.
We also thank Tommaso Toffoli,
Patrizia Mentrasti,
Luca Aceto, Anna Ing\'{o}lfsd\'{o}ttir,
Anders Claesson, and
Magn\'{u}s M\'{a}r Halld\'{o}rsson
for the many helpful suggestions and encouragements.

\end{document}